\documentclass[11pt]{article}
\usepackage{latexsym}
\usepackage{amsfonts}
\usepackage{graphicx}
\usepackage{epstopdf}
\usepackage{xcolor}
\usepackage[sort]{natbib}
\usepackage{lineno}
\newcommand{\coln}{\hspace*{-6pt}{\bf :}}
\newtheorem{theorem}{Theorem}
\newtheorem{lemma}{Lemma}
\newtheorem{corol}{Corollary}
\newtheorem{property}{Property}

\newcommand{\Theorem}[1]{\begin{theorem}\coln ~#1 \end{theorem}}
\newcommand{\Lemma}[1]{\begin{lemma}\coln ~#1 \end{lemma}}

\newcommand{\proof}[1]{\noindent{\bf Proof:~}#1\hfill ~$\Box$}
\newcommand{\be}{\begin{equation}}
\newcommand{\ee}{\end{equation}}

\usepackage[normalem]{ulem}

\begin{document}
%	\linenumbers
	\title{\bf Obnoxious Facility Location: The Case of Weighted Demand Points}
	\author{Pawel Kalczynski\\
		Steven G. Mihaylo College of Business and Economics\\
		California State University-Fullerton\\
		Fullerton, CA 92834.\\e-mail: pkalczynski@fullerton.edu\and Atsuo Suzuki \\
		Department of Mathematical Sciences \\
		Nanzan University \\
		27 Seirei, Seto, Aichi, 489-0863, Japan.\and
   Zvi Drezner\\
		Steven G. Mihaylo College of Business and Economics\\
		California State University-Fullerton\\
		Fullerton, CA 92834.\\e-mail: zdrezner@fullerton.edu
	}
%	\author{}
	\date{}
	\maketitle
	
	\begin{abstract}
The problem considered in this paper is the weighted obnoxious facility location in the convex hull of demand points. The objective function is to maximize the smallest weighted distance between a facility and  a set of demand points. Three new optimal solution approaches are proposed. Two variants of the ``Big Triangle Small Triangle" global optimization method, and a procedure based on intersection points between  Apollonius circles. We also compared the results with a multi-start approach using the non-linear multi-purpose software SNOPT.  Problems with 1,000 demand points are optimally solved in a fraction of a second of computer time.
\end{abstract}

\renewcommand{\baselinestretch}{1.6}
\renewcommand{\arraystretch}{0.625}
\large
\normalsize

\section{Introduction}

The most common objective of obnoxious facilities models \citep[for example, ][]{DKS18,Erk:89} is to maximize the shortest distance to the closest facility. 
In this paper we extend the single obnoxious facility model by introducing weights to demand points. Some demand points are affected differently by the same facility. For example, residential neighborhoods, schools, hospitals, are affected more by noise than commercial neighborhoods. We therefore apply the objective of maximizing the shortest {\it weighted} distance between the facility and demand points rather than just the distance.  One approach to solve the weighted case is applying the weighted Voronoi diagram \citep{AE84} which is a complicated procedure. In this paper we optimally solve the single facility case by relatively simple and very effective approaches.

\begin{figure}[ht]
	\begin{center}\setlength{\unitlength}{1in}
		\begin{picture}(4,4)
		\includegraphics[width=4in]{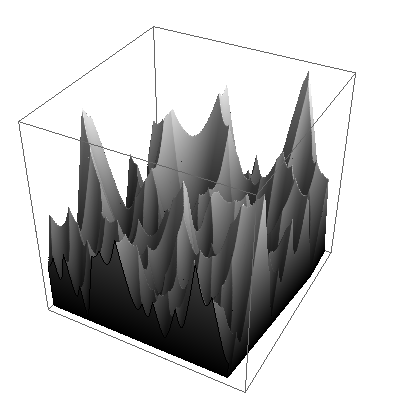}
		\end{picture}
		\caption{\label{surface} The surface of the objective function in the convex hull}
	\end{center}
\end{figure}

The surface of the single weighted obnoxious facility objective function has many local optima. For example, the surface in the convex hull of an $n=100$ demand points instance, tested in the computational experiments, is depicted in Figure \ref{surface}.
As is demonstrated in the computational experiments, non-linear multi-purpose solvers do not perform well and converge to a local optimum depending on the starting solution. Therefore, solution procedures tailored to the special structure of the problem are required for obtaining the optimal solution which is the ``tallest" peak.

Early single  obnoxious facility location models are \citet{CG78,Sham,GD75a}. \citet{Sham} solved a different problem. They found the largest circle which is empty of points. The circle's center is the location of an obnoxious facility that maximizes the shortest distance to a set of points. \citet{Sham} suggested to solve the problem by applying Voronoi diagrams \citep{AKL13,SOK95,OKSC00,V08}. \citet{CG78} solved the problem in a network environment. Later on, \citet{HPT81} assumed that the nuisance generated by the demand points declines by the square of the distance.
They defined the objective of finding the best location for a facility, minimizing the weighted sum of $\frac{1}{d^2}$.

Recent papers on the multiple obnoxious facilities location problem are \citet{DKS18,DDK20}. \citet{DKS18} heuristically solved, the unweighted version of the problem by the standard (unweighted) Voronoi diagram. \citet{DDK20} defined the cooperative obnoxious multi-facilities problem. As in \citet{HPT81}, it is assumed that nuisance generated by the facilities declines by the square of the distance. Each demand point is affected by the total nuisance emitted by the facilities and the objective is to minimize the total nuisance perceived by the most affected demand point.

\section{Heuristic Approach}

To solve the problem by a non-linear solver such as SNOPT \citep{SNOPT}, we formulate it as a non-linear optimization problem. The problem is not convex and many local optima may exist. Let $d_i(X)$ be the distance between the facility located at $X$ and demand point~$i$, and $L$ be the minimum weighted distance to be maximized. The formulation is:
\begin{eqnarray}
&&\max \{~L~\}\nonumber\\
\mbox{Subject to:}\label{form}\\
&&w_id_i(X)\ge L\mbox{ for }i=1,2,\dots,n\nonumber\\
&&A X \le b,\nonumber
\end{eqnarray}
where $A$ and $b$ are the matrix and vector representing the linear constraints of the convex hull of the demand points. We also provided the analytical gradient of the non-linear constraints which significantly improved the performance of SNOPT. 
We applied SNOPT by a multi-start random starting solutions
Unlike the BTST and Apollonius methods described below, this approach does not guarantee finding the global optimum due to the non-convexity of the problem. 

\section{Optimal Algorithms}

We first detail properties that are applied by the optimal algorithms.

\subsection{Finding the optimum for $n=3$ points}

\begin{figure}[htp]
	\begin{center}\setlength{\unitlength}{1in}
		\begin{picture}(2,1.8)
		\includegraphics[width=2in]{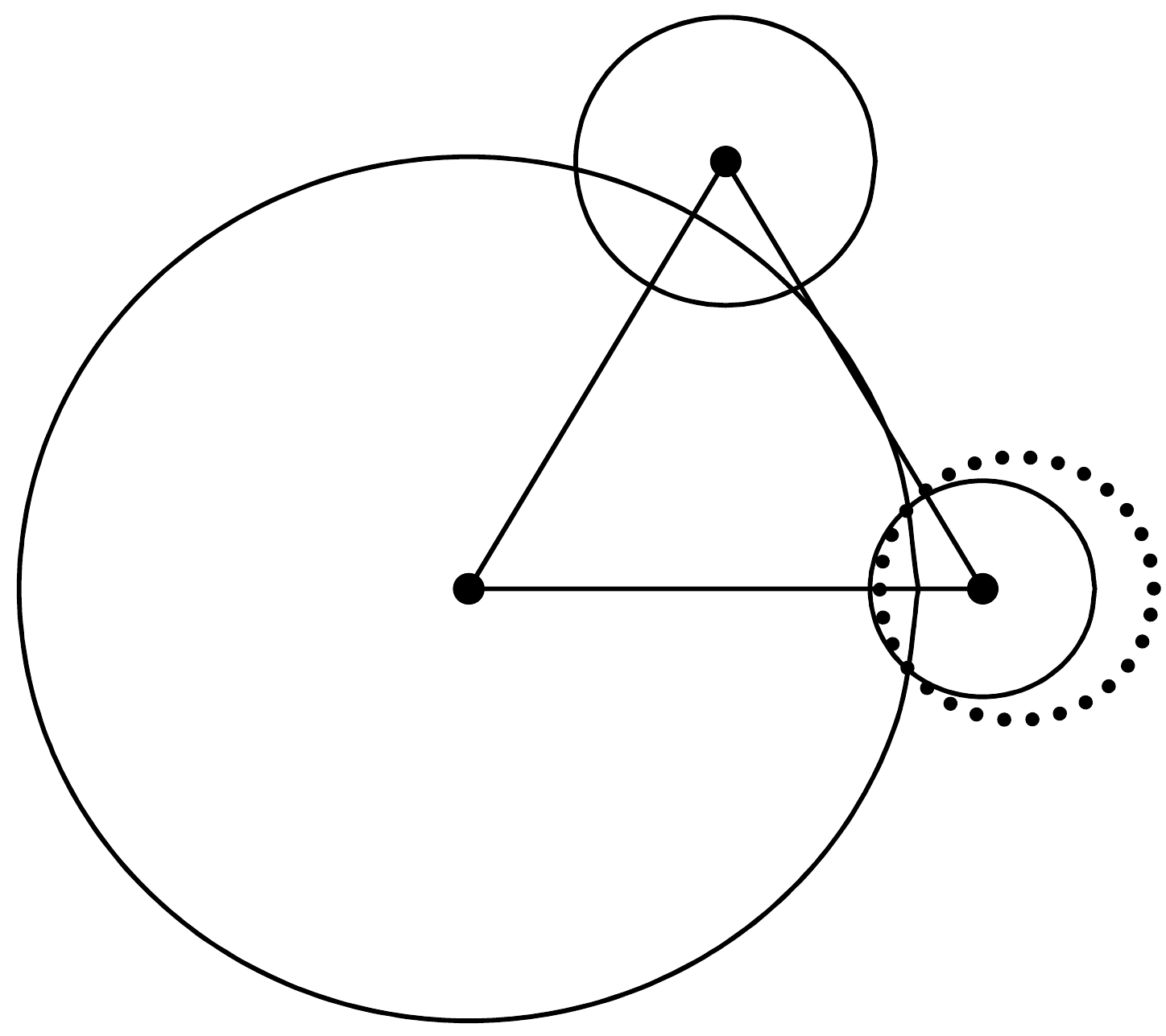}
		\put(-1.32,0.7){\scriptsize{$A$}}
		\put(-0.28,0.7){\scriptsize{$C$}}
		\put(-0.71,1.47){\scriptsize{$B$}}
		\end{picture}
		\caption{\label{circles} Multiple local minima example}
	\end{center}
\end{figure}

To illustrate the problem, consider the following process. Let $L$ be the minimum weighted distance to the vertices of the triangle.  When $L=0$, the feasible region is the whole triangle. As $L$ increases, the feasible region shrinks. It is outside three circles of radii $\frac{L}{w_i}$ centered at the three vertices (see Figure \ref{circles}). The optimal solution is achieved when the feasible region reduces to a point. The problem is not convex and therefore there may be several local minima. Consider Figure~\ref{circles} of an equilateral triangle with vertices A, B, C with corresponding weights of 1, 3, and 4.
One local minimum is on the right side close to point B, and when $L$ is further increased, the feasible area near the top circle disappears, and the optimum is on the right side near point C. The weighted distances between the optimal point and vertices A, C are the same. The optimum is the intersection point between the right side and the Apollonius circle (marked as a dotted circle) based on the vertices A and C.

\begin{figure}[htp]
	\begin{center}\setlength{\unitlength}{1in}
		\begin{picture}(2,1)
		\includegraphics[width=2in]{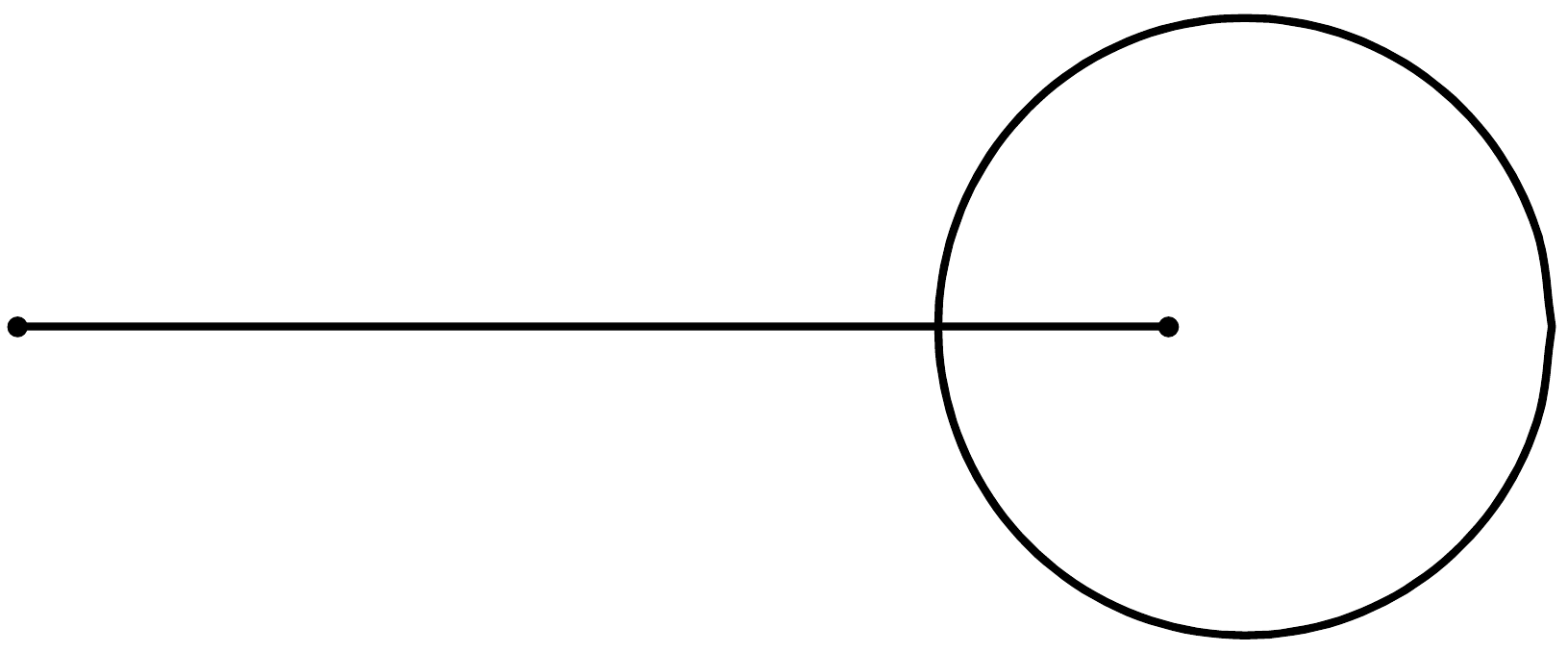}
		\put(-1.65,0.2){Minimum}
		\put(-0.98,0.23){\vector(1,1){0.15}}
		\put(0.15,0.2){Maximum}
		\put(0.15,0.23){\vector(-1,1){0.15}}
		\end{picture}
		\vspace{-0.3in}
		\caption{\label{apol} The Apollonius circle}
	\end{center}
\end{figure}
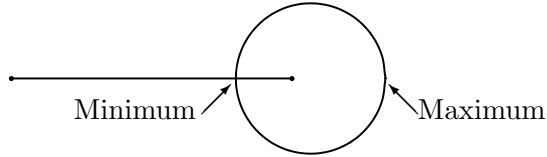

The Apollonius circle \citep{P08,J13} is the set of all points which have the same weighted distance to two points (see Figure \ref{apol}). It is the same condition as ``the ratio between the distances is a constant" (ratio of $\frac{w_1}{w_2}$ or $\frac{w_2}{w_1}$).
All the points that have equal weighted distance from two points $(x_1,y_1)$ and $(x_2,y_2)$ satisfy:
\begin{equation}
w_1^2\left((x-x_1)^2+(y-y_1)^2\right)=w_2^2\left((x-x_2)^2+(y-y_2)^2\right)
\end{equation}
yielding the Apollonius circle:
\begin{equation}\label{eq1}
(w_1^2-w_2^2)(x^2+y^2)-2x\left\{w_1^2x_1-w_2^2x_2\right\}-2y\left\{w_1^2y_1-w_2^2y_2\right\}+w_1^2(x_1^2+y_1^2)-w_2^2(x_2^2+y_2^2)~=~0
\end{equation}
It is the Apollonius circle if $w_1\ne w_2$, and the perpendicular bisector (a circle centered at infinity with an infinite radius) when $w_1=w_2$.

The problem is defined by three points on the vertices of a triangle $T$. A solution $(x,y)$ in the triangle is sought.  The formulation of the $n=3$ problem is similar to the general formulation (\ref{form}):
\begin{eqnarray}
&&\max \{~L~\}\nonumber\\
\mbox{Subject to:}\label{form6}\\
&&w_id_i(X)\ge L\mbox{ for }i=1,2,3\nonumber\\
&&X\in T\nonumber
\end{eqnarray}

\begin{figure}[htp]
	\begin{center}\setlength{\unitlength}{1in}
		\begin{picture}(2,1.5)
		\includegraphics[width=2in]{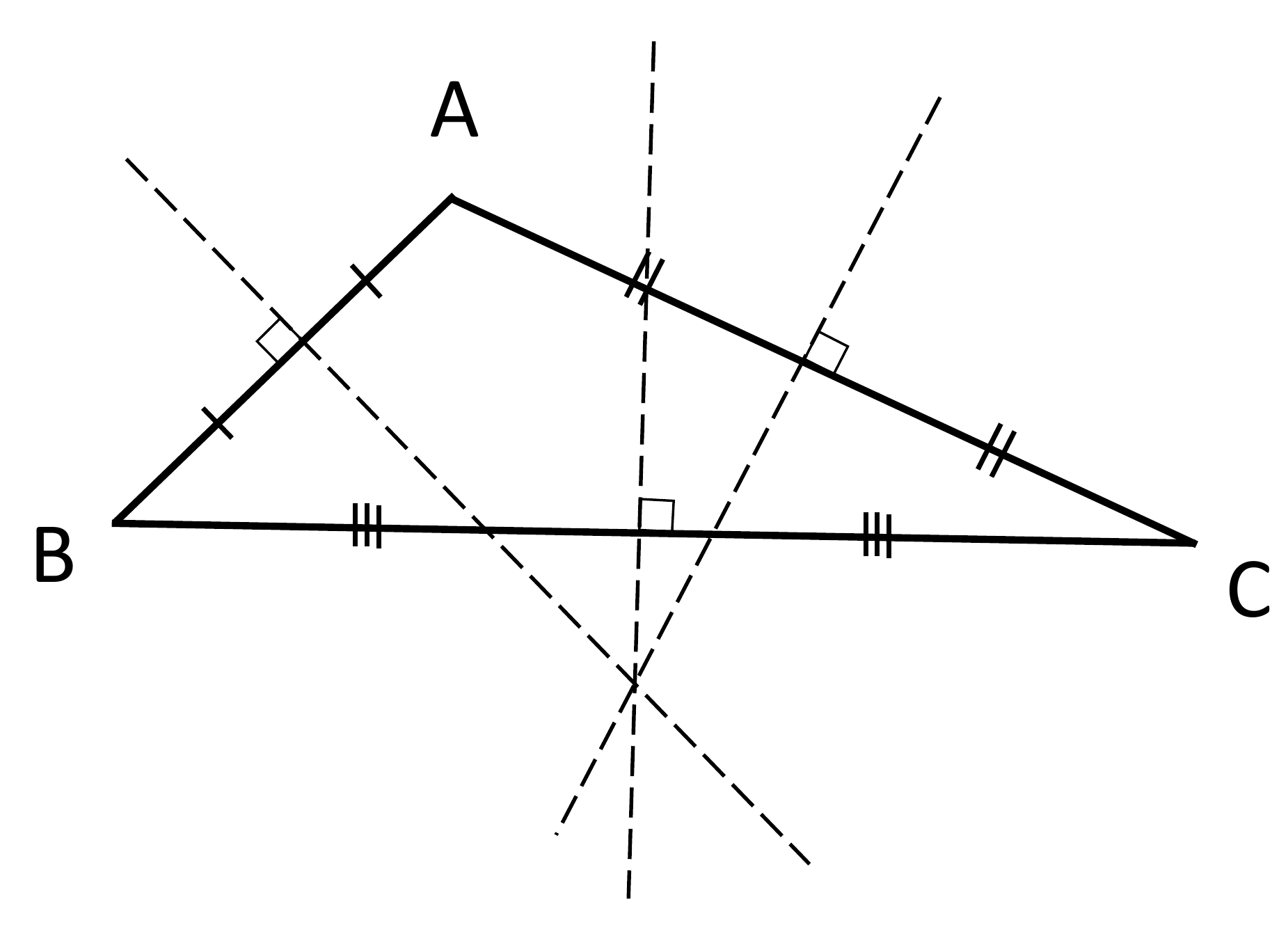}
\put(-0.67,0.35){\scriptsize Optimum}
\put(-0.69,0.4){\vector(-1,1){0.2}}		
		\end{picture}
		\caption{\label{unweighted} An unweighted case}
	\end{center}
\end{figure}

The formulation has six constraints. Of the last three constraints ($X\in T$), a maximum of one constraint can be tight. If two constraints are tight, $X$ is on a vertex that cannot be an optimal solution.
Even when the three weights are equal (the ``standard" unweighted problem), the optimum may be on the side of the triangle. See Figure \ref{unweighted}. The point where all distances are equal is the intersection between the perpendicular bisectors to the sides of the triangle. This intersection point is outside the triangle from below. The optimum is on a side of the triangle at an intersection point between a perpendicular bisector and one of the other sides as marked in the figure.

We first find in each triangle, the point that has equal weighted distance to the three vertices. If all the weights are equal, it is a Voronoi point.
Let the three vertices be $(x_k,y_k)$ and the three weights be $w_k$   for $k=1,2,3$.

The center of the Apollonius circle for $k=1,2$ by equation (\ref{eq1}) is at
\begin{equation}\label{center}
\left(\frac{w_1^2x_1-w_2^2x_2}{w_1^2-w_2^2},\frac{w_1^2y_1-w_2^2y_2}{w_1^2-w_2^2}\right)
\end{equation}
and its radius $R$ is:
\begin{eqnarray}
R&=&\sqrt{\left(\frac{w_1^2x_1-w_2^2x_2}{w_1^2-w_2^2}-\frac{w_1x_1+w_2x_2}{w_1+w_2}\right)^2+\left(\frac{w_1^2y_1-w_2^2y_2}{w_1^2-w_2^2}-\frac{w_1y_1+w_2y_2}{w_1+w_2}\right)^2}\nonumber\\&=&\frac{w_1w_2}{|w_1^2-w_2^2|}\sqrt{(x_1-x_2)^2+(y_1-y_2)^2}\label{radius}
\end{eqnarray}

There are three Apollonius circles, one circle for each of the three pairs of indices.
An intersection point between the circles, if there is one, is at equal weighted distance from the three vertices.
 
\begin{figure}[htp]
 	\begin{center}\setlength{\unitlength}{1in}
 		\begin{picture}(5.4,2.5)
 		\includegraphics[width=2in]{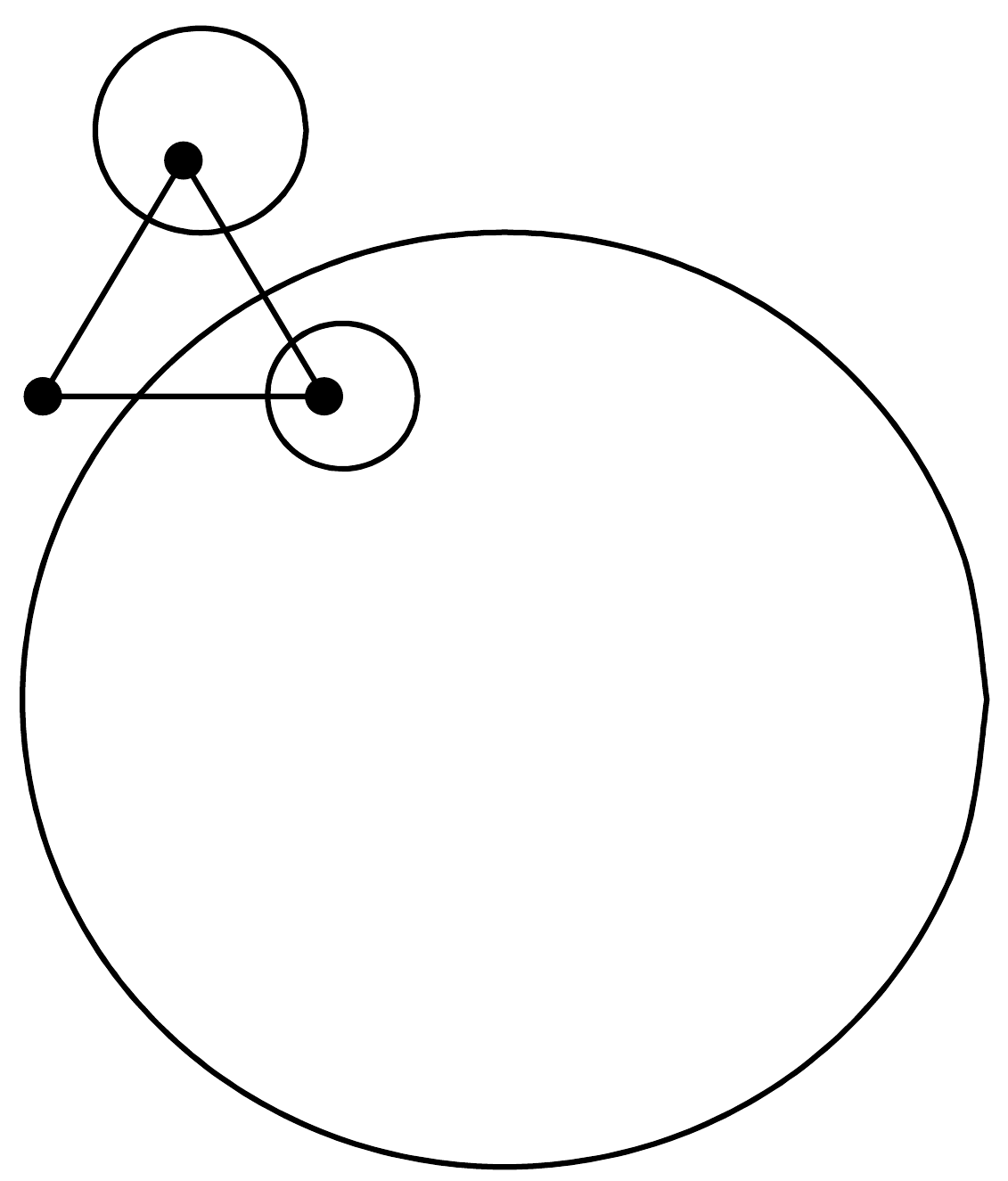}
 		\includegraphics[width=1.8in,height=2.4in]{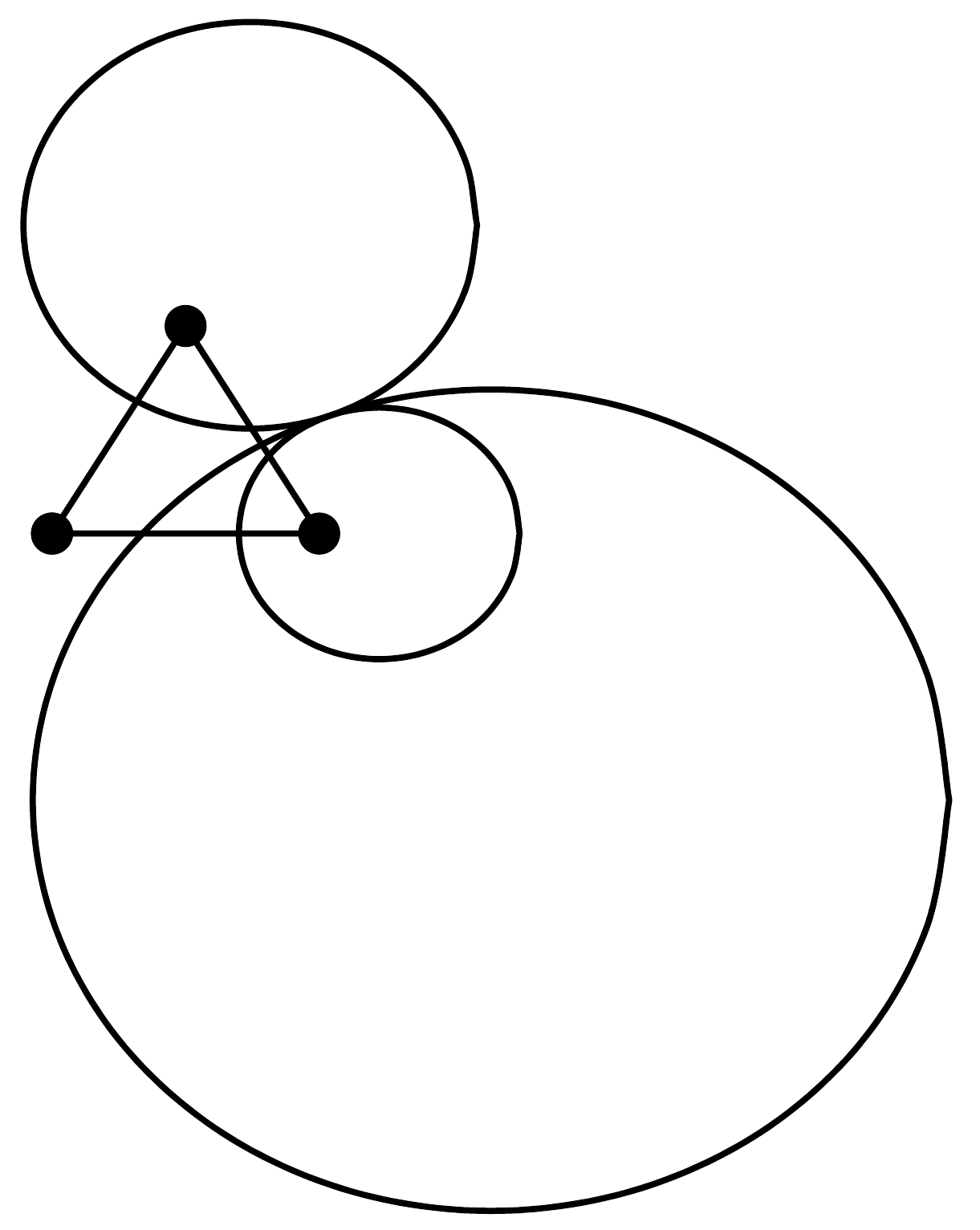}
 		\includegraphics[width=1.8in]{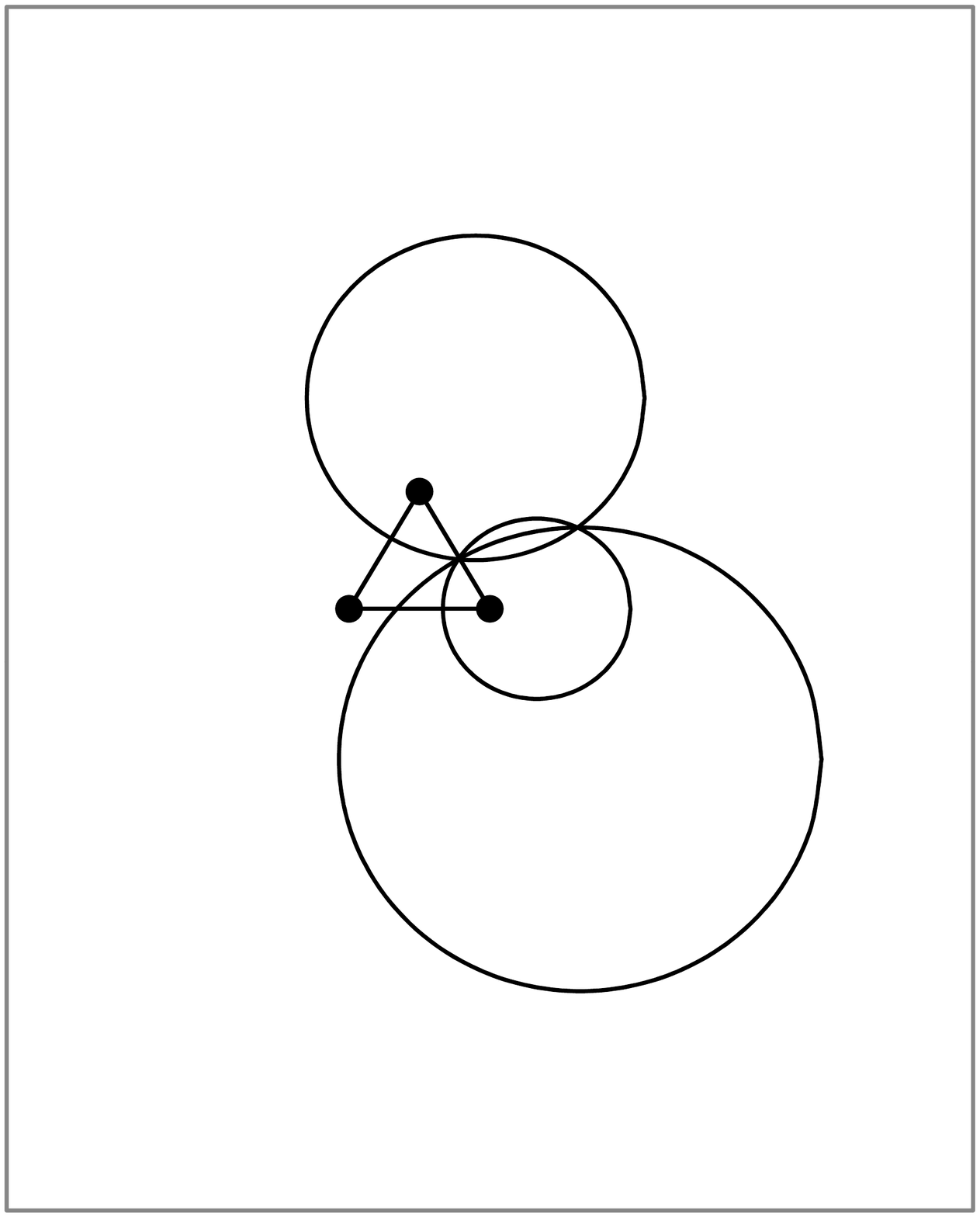}
 		\put(-3.4,-0.15){Weights=$\frac{12}{7}$, 3, 4}
		\put(-1.4,-0.15){Weights=2, 3, 4}
 		\put(-5.4,-0.15){Weights=1, 3, 4}
 		\end{picture}
 		\caption{\label{apols} Apollonius circles}
 	\end{center}
 \end{figure}

See Figure \ref{apols} where the Apollonius circles are drawn for the three vertices depicted in Figure~\ref{circles}. On the left side is the original problem where the weights are 1,3,4. The Apollonius circles do not intersect. We calculated that if we change $w_1$ to $\frac{12}{7}$ and keep the other two weights, the three circles are tangent to one another. This case is depicted at the center of the figure. On the right side, the Apollonius circles are drawn for the case that the first weight is changed from 1 to 2. 
There are two intersection points, one inside the triangle and one outside. Note that if two Apollonius circles intersect, the third one must pass through the intersection points as well because all three weighted distances are the same at an intersection point of any two Apollonius circles. Therefore, it is sufficient to find the intersection point(s) of one pair of Apollonius circles, and find out if there are any. There is no need to find the intersection point(s) between all possible pairs of the Apollonius circles.

\Lemma{\label{outside} Suppose that $w_1\neq w_2$, then the center of the Apollonius circle is on the line connecting the two points outside the segment connecting the two points on the side of the larger weight.}
\proof{By rotation and translation we can get $y_1=y_2=x_1=0$. The center by (\ref{center}) is at $(\frac{w_2^2}{w_2^2-w_1^2}x_2, 0)$ which is on the $x$-axis at a point right to $(x_2, 0)$ when $w_1<w_2$.
}

\Theorem{\label{in} There can be at most one intersection point in the interior of the triangle.}
\proof{Suppose that the smallest weight is $w_1$. Consider Figure \ref{apols} where three cases are plotted. The smallest weight is on the bottom left vertex. Consider the two Apollonius circles based on vertex 1 (bottom left). They are the top two smaller
	 circles. By Lemma \ref{outside}, the centers of the two small circles are on the extensions of the bottom left side of the triangle and the extension of the left side both to the right of the triangle. Therefore, the line connecting them is right to the right side and is outside the triangle. If there are two intersection points (such as the case on the right), the line connecting the two intersection points is perpendicular to the line connecting the two centers. The center of the segment connecting the two intersection points is therefore outside the triangle. If the two points were inside the triangle, the center, as a linear combination of the two points, must be inside the triangle because the triangle is convex. This proves the theorem by contradiction.	
}

The following Lemma is obvious by Figure \ref{apol}.
\Lemma{\label{lem1} Points on the Apollonius circle based on two points have the minimum weighted distance at the intersection point between the Apollonius circle and the segment connecting the two points, and is monotonically increasing along the circle obtaining its maximum weighted distance on the intersection point with the line extension of the segment.}

\Theorem{\label{inside}A local optimum in the interior of the triangle, must have equal weighted distances to all three vertices.}
\proof{If all three weighted distances are different, then moving infinitesimally away from the point with the shortest weighted distance increases the value of the objective function. If the local optimum is at a point at which two weighted distances are equal to $L$ and the weighted distance to the third point is greater than $L$, then moving it infinitesimally along the Apollonius circle increasing $L$ will increase the value of the objective function. By Lemma \ref{lem1}, for a point on the Apollonius circle in the interior of the triangle, there is a direction on the Apollonius circle with increased $L$ because the largest weighted distance on the Apollonius circle is outside the triangle. Therefore, a local optimum inside the triangle must have three equal weighted distances to all three vertices and is an intersection point between any two Apollonius circles.
}

\Theorem{\label{side}A local optimum on a particular side of the triangle can only be at the intersection point between the side and an Apollonius circle based on one of the end vertices of that side and the third vertex.}
\proof{A vertex of the particular side has an objective of $L=0$. As the point moves away from that vertex, $L$ increases until it matches the weighted distance to the third vertex (after which it decreases), which is the intersection point between the Apollonius circle based on this vertex and the third vertex and it can be a local optimum on that side. Another possibility is that it reaches the weighted distance to the second vertex and the weighted distance to the third vertex is greater than $L$. This is not a local optimum because moving the point infinitesimally on the Apollonius circle based on the two vertices of the particular side, increases $L$ in any direction (in particular in a direction into the interior of the triangle), because this point has the smallest weighted distance by Lemma \ref{lem1}.
}

\Theorem{\label{glob}If there is an intersection point between two Apollonius circles in the interior of the triangle, it is the global optimum.}
\proof{By the triangle inequality, the sum of the distances between the intersection point and any two vertices is greater than the side of the triangle connecting the vertices. Since the weighted distances to the two vertices are the same, the smaller weighted distance on the side of the triangle must be smaller. Therefore, this interior point has a greater value of the objective function than any point on a side of the triangle. By Theorem~\ref{in}, there is at most one such intersection point and thus it is the global optimum.
}

\subsubsection{\label{suggest}Suggested Solution Procedure for $n=3$}

Detailed description of all the needed calculations in the suggested procedure are given in the Appendix.

\begin{enumerate}
	\item Find two Apollonius circles by equations (\ref{center}) and (\ref{radius}).
	\item \label{st2}Find the intersection points between the two Apollonius circles if there are any.
	\item If there are such intersection points, check whether an intersection point is inside the triangle. There can be at most one such intersection point by Theorem \ref{in}. This point is the global optimum by Theorem \ref{glob}. Stop.
	\item Otherwise, there is no local optimum in the interior of the triangle.
	\begin{enumerate}
		\item Find the third Apollonius circle.
		\item	 Find the intersection points between each Apollonius circle with the two sides that include the third vertex.
		\item Calculate the value of the objective function at all intersection points with the sides, and the best one is the optimal solution by Theorem \ref{side}.
	\end{enumerate}
\end{enumerate}

\subsection{The General Case}

We distinguish between two cases; (i) the optimal solution is in the interior of the convex hull, or (ii) the optimal solution is on the periphery of the convex hull.

\Lemma{\label{lem3}If the solution is in the interior of the convex hull, it is at an intersection point between two Apollonius circles with a common vertex.}
\proof{There is at least one point $i$ whose weighted distance to the solution point is the maximum possible $L$. If there is no second point with weighted distance $L$, a small movement of the solution point away from point $i$ will increase $L$. Therefore, for an interior optimal solution there must be at lease two points $i,j$ whose weighted distance is $L$. Therefore, the optimal solution is on the Apollonius circle based on $i,j$. If the weighted distance to all other points is greater than $L$, moving the solution point along the Apollonius circle will increase $L$, unless the solution point is at the farthest point of the Apollonius circle (see Figure \ref{apol}). However, in this case moving away from points $i$ and $j$ (because the farthest point is on the extension of the line connecting the points) will increase the distance to both $i$ and $j$. Therefore, there must be a third point $k$ with a weighted distance $L$. In conclusion, an optimal solution in the interior of the convex hull must be at an intersection point between two Apollonius circles with a common vertex.
}
\Lemma{\label{lem4}If the solution is on a side of the convex hull, it must be on an intersection point between the side and an Apollonius circle.}
\proof{The solution point cannot be on an endpoint of the side, because the weighted distance to the endpoint is zero. Suppose that the solution point is at least a weighted distance $L$ from all points (including points which are not the ends of the side). There is at least one point $i$ at weighted distance $L$. If the weighted distance to all other points is greater than $L$, the solution point can be moved on the side away from point $i$, thus increasing $L$. Therefore there must be a second point $j$ at weighted distance $L$. In conclusion, the best solution point on a side of the convex hull must be on an intersection point between the side and an Apollonius circle based on two points $i$ and $j$.
}

These two Lemmas lead to the following theorem.

\Theorem{The solution point to the single obnoxious facility location problem with weighted distances is either on an intersection point of an Apollonius circle and a side of the convex hull, or at an interior point at the intersection of two Apollonius circles with a common vertex.}

\subsection{Applying the Standard Delaunay Triangulation}

\begin{figure}[ht!]
	\begin{center}
		\setlength{\unitlength}{1in}
	\begin{picture}(4,4)
		\includegraphics[width=4in]{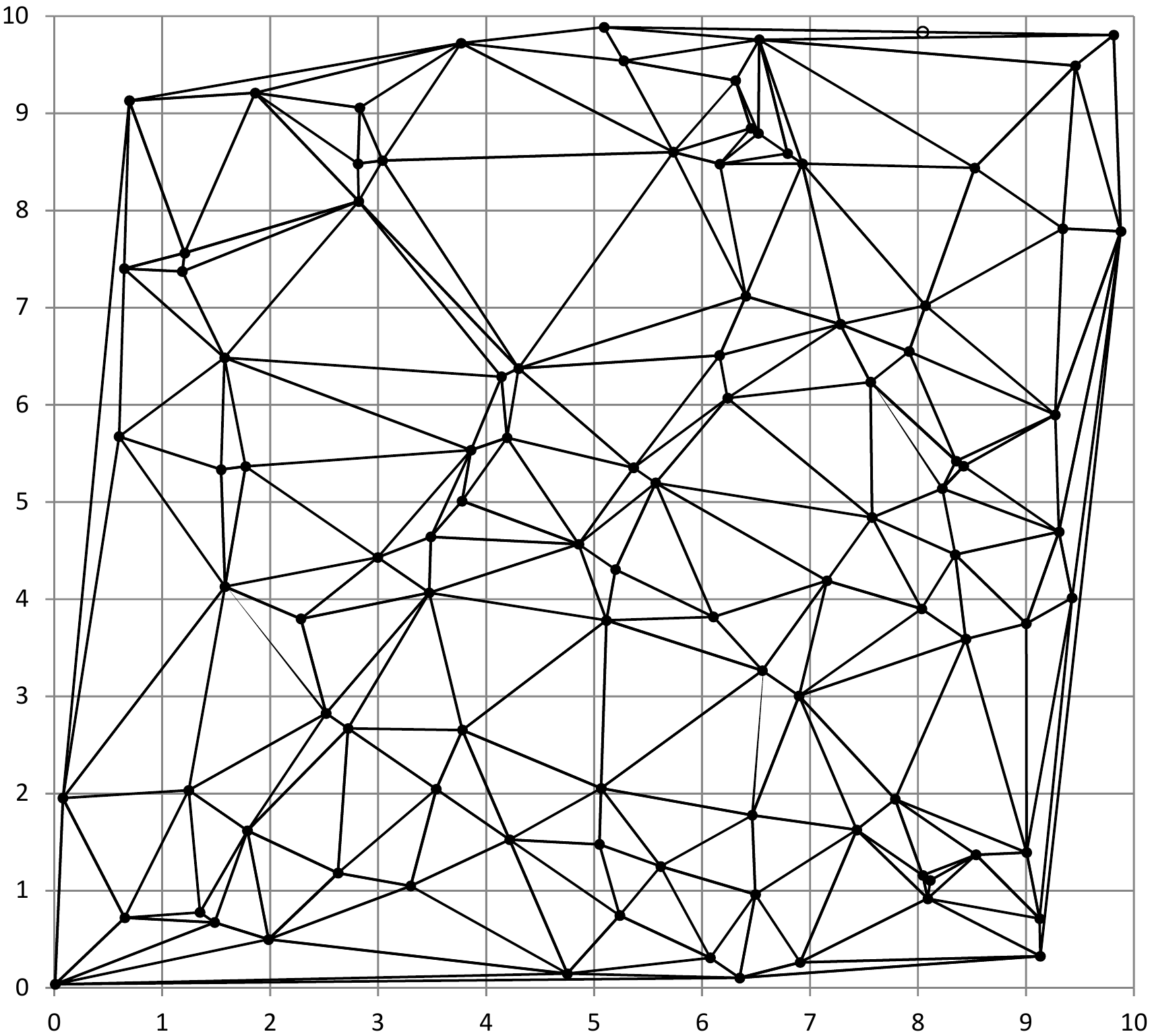}
\put(-0.57,3.57){\tiny Optimum}
\put(-0.58,3.6){\vector(-2,-1){0.2}}		
\end{picture}
\caption{\label{Dela}The Delaunay Triangulation for the $n=100$ instance}
\end{center}
\end{figure}

We triangulate the demand points as vertices by the standard (non-weighted) Delaunay triangulation \citep{LS80,AKL13}. The vertices of each triangle are demand points and thus have associated weights.
% We first show how to find the optimal solution in one triangle considering only its three vertices in the objective function.
See Figure \ref{Dela} for the Delaunay triangles for the $n=100$ tested instance. The small empty circle near the top of the figure on the right is the optimal solution.

\subsubsection{\label{BTST}The BTST Based Procedure}
A relative accuracy $\epsilon>0$ is given. The BTST (Big Triangle Small Triangle) global optimization procedure \citep{DS04} is processed as follows. 
\begin{itemize}
\item A list of triangles \citep[Delaunay triangles, ][]{AKL13,OKSC00} covering the feasible region is created. There are $V$ triangles for $i=1,\ldots,V$.
\item A lower bound, for each triangle $1\le i\le V$, $LB_i$ and an upper bound $UB_i$ are found. The largest lower bound is $\overline{LB}$.
\item All triangles for which  $UB_i\le \overline{LB}(1+\epsilon)$ are discarded.
\item The triangle with the largest $UB_i$ is selected and is split into four ``small" triangles.
\item The lower and upper bounds are found in each ``small" triangle.
\item $\overline{LB}$ may be updated.
\item  All ``small" triangles for which $UB_i\le \overline{LB}(1+\epsilon)$, and the ``big" triangle, are discarded.
\item The process continues until all triangles are eliminated.
\end{itemize}

\subsubsection{First Set of Bounds}
Consider a triangle whose vertices are not necessarily demand points, and thus have no associated weights.

\noindent{\bf A lower bound $LB^{(1)}$} is the value of the objective function at any point in the triangle, such as its (unweighted) center of gravity.

\noindent{\bf An upper bound $UB^{(1)}$} is found as follows. For each demand point $1\le i\le n$, find the largest distance to the three vertices of the triangle (defined as $d_i$). The minimum $w_id_i$ for $1\le i\le n$ is the upper bound $UB^{(1)}$.

\subsubsection{Second Set of Bounds}

The second set of bounds can be used only if all the vertices of the triangle are demand points. The optimal solution to the $n=3$ vertices problem is found.

\noindent{\bf A lower bound $LB^{(2)}$} is the value of the objective function based on all $n$ demand points at the solution point to the $n=3$ problem.

\noindent{\bf An upper bound $UB^{(2)}$} is the value of the objective function based on $n=3$ demand points at the solution point to the $n=3$ problem.

\subsubsection{Specific BTST procedures}

\subsubsection*{BTST1}

For this approach there is no need to find the optimal solution for the $n=3$ problems. We apply the first set of bounds. $LB^{(1)}$, $UB^{(1)}$ in the BTST procedure.

\subsection*{BTST2}

Different lower and upper bounds are applied in the initial phase and in subsequent generation of ``small" triangles.

\begin{enumerate}
	\item Create the list of $t$ {\it standard} Delaunay triangles \citep{AKL13,OKSC00}. By design, the three vertices of each triangle are demand points and the union of all triangles is the convex hull of the demand points.
	\item For each triangle $1\le i\le t$ find the lower and upper bounds $LB_i$ and $UB_i$ by  the second set of bounds $LB^{(2)}$, $UB^{(2)}$.
	\item Let $\overline{LB}$ be the maximum among all $LB_i$.
	\item Discard from the list all the triangles for which $UB_i\le \overline{LB}(1+\epsilon)$.
	\item Apply BTST on the remaining list of triangles using the first set of bounds.
Note that it is possible that the list of triangles is empty if in the triangle with the best lower bound,  all weighted distances to demand points outside the specific triangle are greater than the optimal solution based on the three vertices of that triangle.

\end{enumerate}

\subsection{The Apollonius Method}

By Lemmas \ref{lem3} and \ref{lem4}, we designed the following algorithm to find the global optimum, i.e., the location for a single facility in the convex hull with the largest weighted minimum distance to the demand points.

\begin{enumerate}

\item Find the sides of the convex hull polygon.

\item Find the Apollonius circles for all pairs of weighted demand points. Find the intersection points between these circles and the sides of the convex hull polygon. Calculate the value of the objective function at each intersection point.

\item For each triangle (triplet) of weighted demand points, find the intersection points of Apollonius circles, if they intersect. Calculate the value of the objective function at intersection points that are inside the triangle. By Theorem \ref{in} there can be at most one intersection point in the interior of the triangle.

\item The solution point is the intersection point with the largest value of the objective function (minimum weighted distance to the demand points).

\end{enumerate}

This approach is guaranteed to find the global optimal solution. Its worst-case complexity is $O(n^4)$.  However, a significant speed-up is achieved in the last step if the current-best value of the objective function is kept in memory. During the procedure of evaluating the value of the objective function, if the minimum weighted distance between the analyzed candidate point and a demand point is less than the best found value, there is no need to evaluate more distances, and the intersection point is skipped and not considered for the optimum. A regression analysis of the run times of the test instances reported in the computational experiments  section found a complexity of $O(n^{2.85})$ with a $p$-value of $3.7\times 10^{-10}$ which is highly significant.

\subsection{Weighted Voronoi Diagrams}

\citet{AE84}  proposed a procedure to optimally generate the weighted Voronoi diagram. Their diagram is not restricted to the convex hull and therefore we need to check whether each candidate point is in the convex hull and calculate the intersections of the Apollonius circles with the boundary of the convex hull. They prove that their algorithm is of complexity $O(n^2)$ and do not report computational experiments. However, checking whether the Voronoi points are in the convex hull, and evaluating all the intersection points between the Apollonius circles with the sides of convex hull,  requires more than $O(n^2)$ time. The objective function needs to be calculated for every intersection point and requires $O(n)$ time. There are $O(n^2)$ Apollonius circles. Each Apollonius circle may intersect with each of the sides of the convex hull. The number of sides of the convex hull increases as $n$ increases.  

Since our proposed optimal algorithms are not complicated to apply, are of polynomial complexity, and performed very well in the computational experiments (optimally solving the $n=1,000$ instance took 0.12 seconds), we opted not to attempt to program and test \citet{AE84}'s  complicated algorithm.
\section{Computational Experiments}

We generated problems with $n=100,200,\ldots,1000$ demand points. In order to allow for future comparisons, the problems were generated by the method in \cite{DKS18}. A sequence $r_k$ of integer numbers in the open range (0, 100,000) is generated. A starting seed $r_1$, which is the first number in the sequence, and a multiplier $\lambda$ which is an odd number not divisible by 5, are selected. We used $\lambda=12219$. The sequence is generated by the following rule for $k\ge 1$:

$$
r_{k+1}=\lambda r_k-\lfloor\frac{\lambda r_k}{100000}\rfloor\times 100000.
$$
The random number between 0 and 10 is $\frac{r_k}{10,000}$.

For demand points (with coordinates between 0 and 10) the $x$ coordinates were generated by $r_1=97$, and for the $y$-coordinates we used $r_1=367$. For the weights we used $r_1=12347$, and $w_i=1+\frac{r_i}{100,000}$ so $1<w_i<2$.
For illustration, the first 30 points are depicted in Table \ref{points}.

\begin{table}[ht!]
	\begin{center}
		\caption{\label{points}The First 30 Randomly Generated Points}
		\medskip
		%	\small
		\setlength{\tabcolsep}{3.5pt}
		\begin{tabular}{|c|c|c|c||c|c|c|c||c|c|c|c|}
			\hline
			$i$&$x$&$y$&$w$&$i$&$x$&$y$&$w$&$i$&$x$&$y$&$w$\\
			\hline
1&	0.0097&	0.0367&	1.12347&	11&	9.2697&	5.8967&	1.54947&	21&	6.5297&	9.7567&	1.77547\\
2&	8.5243&	8.4373&	1.67993&	12&	6.4643&	1.7773&	1.97393&	22&	6.4043&	7.1173&	1.46793\\
3&	8.4217&	5.3687&	1.06467&	13&	7.2817&	6.8287&	1.45067&	23&	4.1417&	6.2887&	1.63667\\
4&	4.7523&	0.1453&	1.20273&	14&	5.0923&	9.8853&	1.73673&	24&	7.4323&	1.6253&	1.47073\\
5&	8.3537&	5.4207&	1.15787&	15&	2.8137&	8.4807&	1.10387&	25&	5.2737&	9.5407&	1.84987\\
6&	3.8603&	5.5333&	1.01353&	16&	0.6003&	5.6733&	1.18753&	26&	9.3403&	7.8133&	1.56153\\
7&	9.0057&	1.3927&	1.32307&	17&	5.0657&	2.0527&	1.42907&	27&	9.1257&	0.7127&	1.33507\\
8&	0.6483&	7.4013&	1.59233&	18&	7.7883&	1.9413&	1.80633&	28&	6.9283&	8.4813&	1.22033\\
9&	1.5777&	6.4847&	1.68027&	19&	5.2377&	0.7447&	1.54627&	29&	6.8977&	3.0047&	1.21227\\
10&	7.9163&	6.5493&	1.21913&	20&	9.4563&	9.4893&	1.87313&	30&	2.9963&	4.4293&	1.72713\\
\hline
		\end{tabular}
	\end{center}
\end{table}

The  optimal algorithms were coded in FORTRAN and used double precision arithmetic. The programs were compiled by an Intel 11.1 FORTRAN Compiler with no parallel processing. They were run on a desktop with the Intel i7-6700 3.4GHz CPU processor and 16GB RAM. Only one processor was used.
The multi-purpose non-linear solver SNOPT \citep{SNOPT} was run on an virtualized Windows server with 16 CPUS and 128 GB of RAM.

In Table \ref{opt} the coordinates of the optimal solutions and the optimal value of the objective function are given. We also depict the number of Delaunay triangles, and the number of sides in the convex hull of the demand points.

\begin{table}[ht!]
	\begin{center}
		\caption{\label{opt}Properties and optimal solutions of the tested instances}
		\medskip
		%	\small
		\setlength{\tabcolsep}{4pt}
		\begin{tabular}{|c|c|c|c|c|c||c|c|c|c|c|c|}
			\hline
			$n$&$x^*$&$y^*$&objective&$\dagger$&$\ddagger$&$n$&$x^*$&$y^*$&objective&$\dagger$&$\ddagger$\\
			\hline
100&	8.04233&	9.83530&	2.13972&	189&	9&	600&	0.04420&	7.14163&	1.04703&	1183&	15\\
200&	9.89778&	3.12986&	1.63585&	387&	11&	700&	0.04420&	7.14163&	1.04703&	1382&	16\\
300&	4.11567&	7.65730&	1.37183&	585&	13&	800&	0.04420&	7.14163&	1.04703&	1582&	16\\
400&	8.88491&	9.85960&	1.10596&	786&	12&	900&	0.04420&	7.14163&	1.04703&	1781&	17\\
500&	0.04420&	7.14163&	1.04703&	985&	13&	1000&	0.04421&	7.14310&	1.04609&	1981&	17\\
			\hline
			\multicolumn{12}{l}{$\dagger$ Number of Delaunay triangles.}\\
			\multicolumn{12}{l}{$\ddagger$ Number of sides in the convex hull.}\\
		\end{tabular}
	\end{center}
\end{table}

\subsection{Solving by BTST}

The performance of the BTST algorithms is depicted in Table \ref{sol}. A relative accuracy of $\epsilon=10^{-10}$ was used. For each BTST variant we show
\begin{itemize}
\item The lower and upper bounds at the end of the first phase,  which scans all the Delaunay triangles.
\item The number of triangles remaining at the end of the first phase.
\item The maximum number of triangles during the second phase of the algorithm (scanning these remaining triangles).
\item  The number of iterations executed at the second phase.
\item The time in seconds to solve the problem by the FORTRAN program.
\end{itemize}

The second set of bounds performs better than the first  set. However, since almost the whole time is consumed by the second phase, which uses the first set of bounds in both variants, the better performance is hardly reflected in the run times. 

\begin{table}[ht!]
	\begin{center}
		\caption{\label{sol}Performance by the BTST methods}
		\medskip
		%	\small
		\setlength{\tabcolsep}{3pt}
		\begin{tabular}{|c|c||c|c|c|c|c|c||c|c|c|c|c|c|}
			\hline
			$n$&Optimal&\multicolumn{6}{c||}{BTST1}&\multicolumn{6}{c|}{BTST2}\\
			\cline{3-14}
			&objective&(1)&(2)&(3)&(4)&(5)&(6)&(1)&(2)&(3)&(4)&(5)&(6)\\
			\hline
100&	2.13972&	1.65610&	4.84069&	84&	101&	1,087&	0.002&	2.00742&	3.63158&	11&	88&	1,044&	0.002\\
200&	1.63585&	1.54729&	5.52649&	72&	510&	7,166&	0.022&	1.61199&	5.12565&	13&	510&	7,086&	0.022\\
300&	1.37183&	1.14756&	5.01672&	155&	219&	388&	0.002&	1.32627&	5.29483&	17&	108&	320&	0.002\\
400&	1.10596&	0.94573&	3.94609&	254&	373&	543&	0.004&	1.10233&	3.43307&	20&	76&	392&	0.002\\
500&	1.04703&	0.88090&	3.94609&	253&	335&	489&	0.005&	0.90744&	3.61715&	25&	121&	328&	0.003\\
600&	1.04703&	0.93155&	3.94609&	177&	285&	1,409&	0.013&	0.95256&	3.61715&	27&	119&	1,294&	0.011\\
700&	1.04703&	0.93155&	3.94609&	163&	264&	1,360&	0.015&	0.95256&	3.61715&	25&	115&	1,279&	0.012\\
800&	1.04703&	0.93155&	3.94609&	128&	218&	1,307&	0.018&	0.95256&	3.03222&	25&	106&	1,249&	0.014\\
900&	1.04703&	0.93155&	3.94609&	91&	171&	1,253&	0.020&	0.95256&	3.03222&	24&	98&	1,221&	0.016\\
1000&	1.04609&	0.93155&	3.94609&	75&	730&	9,004&	0.127&	0.91536&	2.36876&	27&	730&	8,981&	0.122\\
			\hline
\multicolumn{14}{l}{(1) Lower bound at the end of the first phase.}\\
\multicolumn{14}{l}{(2) Upper bound at the end of the first phase.}\\
\multicolumn{14}{l}{(3) Number of triangles remaining at the end of the first phase.}\\
\multicolumn{14}{l}{(4) Maximum number of triangles throughout the iterations in the second phase.}\\
\multicolumn{14}{l}{(5) Number of iterations in the second phase.}\\
\multicolumn{14}{l}{(6) Run time in seconds.}\\
		\end{tabular}
	\end{center}
\end{table}

\begin{table}[ht!]
	\begin{center}
		\caption{\label{ap}Performance by the Apollonius Method and SNOPT}
		\medskip
		%	\small
		%\setlength{\tabcolsep}{2.5pt}
				\setlength{\tabcolsep}{7pt}	
		\begin{tabular}{|c||c|c|c|c||c|c|}
			\hline
			$n$&\multicolumn{4}{c||}{Apollonius Method}&\multicolumn{2}{|c|}{SNOPT}\\
			\cline{2-7}
			&(1)&(2)&(3)&(4)&(5)&(6)\\
			\hline
100&	44,550&	161,700&	49,951&	0.02&	312&	29.30\\
200&	218,900&	1,313,400&	375,223&	0.08&	201&	43.18\\
300&	583,050&	4,455,100&	1,198,954&	0.27&	4&	53.56\\
400&	957,600&	10,586,800&	2,815,919&	0.70&	148&	67.25\\
500&	1,621,750&	20,708,500&	5,489,533&	1.42&	3&	80.13\\
600&	2,695,500&	35,820,200&	9,376,817&	2.48&	3&	93.84\\
700&	3,914,400&	56,921,900&	14,873,925&	4.00&	3&	109.52\\
800&	5,113,600&	85,013,600&	22,079,425&	5.67&	3&	126.43\\
900&	6,877,350&	121,095,300&	31,276,311&	8.05&	3&	138.58\\
1000&	8,491,500&	166,167,000&	42,809,302&	11.03&	3&	154.69\\
			\hline
			\multicolumn{7}{l}{(1) Number of potential intersection points with the convex hull.}\\
			\multicolumn{7}{l}{(2) Number of triplets of demand points.}\\
			\multicolumn{7}{l}{(3) Number of times the objective function was calculated.}\\
			\multicolumn{7}{l}{(4) Time in seconds.}\\
			\multicolumn{7}{l}{(5) Number of times out of 10,000 that optimum found.}\\
			\multicolumn{7}{l}{(6) Time in seconds for all 10,000 runs.}\\
		\end{tabular}
	\end{center}
\end{table}

\begin{table}[ht!]
	\begin{center}
		\caption{\label{square}Solving the test problems in a square}
		\medskip
		%	\small
		%\setlength{\tabcolsep}{2.5pt}
		\setlength{\tabcolsep}{7pt}	
		\begin{tabular}{|c||c|c|c||c|c||c|c|}
			\hline
			$n$&$x^*$&$y^*$&Objective&\multicolumn{2}{c||}{Apollonius Method}&\multicolumn{2}{|c|}{SNOPT}\\
			\cline{5-8}
			&&&&(1)&(2)&(3)&(4)\\
			\hline
			100&	10.00000&	2.77952&	2.25773&	49,927&	0.02&	230&	34.50\\
			200&	10.00000&	3.09849&	1.69987&	375,219&	0.08&	187&	36.23\\
			300&	10.00000&	2.90534&	1.41960&	1,198,996&	0.25&	143&	49.18\\
			400&	8.88615&	10.00000&	1.14393&	2,816,019&	0.69&	191&	60.19\\
			500&	0.00000&	7.17256&	1.09468&	5,489,777&	1.36&	3&	72.30\\
			600&	0.00000&	7.17256&	1.09468&	9,377,117&	2.33&	3&	83.03\\
			700&	0.00000&	7.17256&	1.09468&	14,874,359&	3.69&	3&	98.98\\
			800&	0.00000&	7.17256&	1.09468&	22,080,139&	5.48&	3&	108.63\\
			900&	0.00000&	7.17256&	1.09468&	31,277,241&	7.80&	3&	123.40\\
			1000&	0.00000&	7.17256&	1.09468&	42,810,610&	10.64&	3&	131.64\\
			\hline
			\multicolumn{8}{l}{(1) Number of times the objective function was calculated.}\\
			\multicolumn{8}{l}{(2) Time in seconds.}\\
			\multicolumn{8}{l}{(3) Number of times out of 10,000 that optimum found.}\\
			\multicolumn{8}{l}{(4) Time in seconds for all 10,000 runs.}\\
		\end{tabular}
	\end{center}
\end{table}

\subsection{Solving by the Apollonius Method}

In Table \ref{ap} the data about the performance of the Apollonius method is depicted:
\begin{itemize}
\item The number of potential intersection points between all Apollonius circles and the side of the convex hull. The actual number of intersection points is much lower. However, we must find out whether each Apollonius circle intersects with every side.
\item The number of triplets of demand points forming triangles ($\frac{n(n-1)(n-2)}{6}$), generating two Apollonius circles each which may intersect.
\item The number of actual calls to the subroutine calculating the value of the objective function for each intersection point. Such calculation requires the evaluation of weighted distances to all demand points. However, if a weighted distance to a demand point is smaller than the best solution found so far, the calculation is terminated because such an intersection point cannot have an improved objective.
\item The run time in seconds.
\end{itemize}

In Table \ref{square} the results of solving the problem in a square of side 10 rather than in the convex hull are reported. The number of sides of the feasible area is 4 rather than the number of sides in the convex hull reported in Table \ref{opt}. The performance of the Apollonius method is about the same as solving the problem in the convex hull of the demand points.

\subsection{Solving by SNOPT}

A multi-start approach from random locations in the square was used to determine the best solution. When SNOPT was re-started from 1,000 randomly generated starting solutions, the optimal solution was missed for several instances. We therefore experimented with 10,000 starting solutions. The performance is depicted in Table \ref{ap}. The optimal solution was found between 3 and 312 times out of 10,000. Run times are longer than required by other approaches and the optimal solution is not guaranteed.

We also report in Table \ref{square} the performance of SNOPT  for the problems in a square rather than the convex hull. SNOPT performed a bit faster because the number of constraints is reduced.

Note that the same 10,000 random starting solutions were applied for each instance. Since the optimal solution is the same for the $n\ge 500$ instances, the number of times the optimal solution found  is also the same for $n\ge 500$.

\section{Conclusions}

We propose optimal solution algorithms for the planar obnoxious facility location problem with weighted distances. The objective is to maximize the minimum weighted distance to all demand points. The location of the facility is restricted to the convex hull of the demand points. Any finite region can also be considered. Otherwise, the solution is at infinity.

The simplest approach to implement is to apply non-linear optimization software, such as SNOPT \citep{SNOPT} in a mult-start approach. However, the optimal solution is not guaranteed and it takes the longest run time.

Two variants of the ``Big Triangle Small Triangle"  \citep[BTST, ][]{DS04} are also quite easy to implement once the software for BTST is available. Only lower and upper bounds in a triangle are required. The bounds for version BTST1   are very easy to implement and the optimal solution (within $\epsilon=10^{-10}$ relative accuracy) is found in 0.13 seconds for $n=1,000$ demand points. The bounds for BTST2 require finding the optimal solution in a triangle The algorithm is more efficient than BTST1 but run time is reduced very little (0.12 seconds for $n=1,000$).

A third approach, termed the Apollonius method, is the most complicated to implement. It finds the intersection points between Apollonius circles and the boundary of the convex hull, and the intersection points between pairs of Apollonius circles. One of these intersection points is the optimal solution. This method takes longer than the BTST variants (about 11 seconds for the $n=1,000$ instance) but the optimal solution is precisely calculated by an explicit formula, and its accuracy is actually dependent on the accuracy of the computer processor.

{\samepage{{\begin{center}{\bf\large Appendix: Detailed Calculations}\end{center}}
\subsection*{Translation and Rotation}}}

We translate and rotate the system of coordinates so that a given point $A_1=(a_1,b_1)$ is transformed to $(0, 0)$ and a given point $A_2=(a_2,b_2)$ is transformed to $(d,0)$ where $d$ is the distance between the two points.
We use matrix notation, i.e., a point $(x,y)$ is denoted as the vector $X=\left(x\atop y\right)$. The distance between the two points $d$ satisfies
\begin{equation}\label{d}
d^2=(a_2-a_1)^2+(b_2-b_1)^2=(A_2-A_1)^T(A_2-A_1).
\end{equation}
A translation and rotation formula of $X$ into $\bar X$ is
\begin{equation}\label{trans}
\bar X=M(X-A_1)
\mbox{ where: }
M=\left(\begin{array}{cc}
\frac{a_2-a_1}{d}&\frac{b_2-b_1}{d}\\
{-\frac{b_2-b_1}{d}}&{\frac{a_2-a_1}{d}}
\end{array}\right) 
\end{equation}

It can be easily verified by using (\ref{d}), that $\bar A_1$ is (0, 0), and $\bar A_2$ is the point $(d,0)$.
\Lemma{\label{MTM}$M^TM=\left(\begin{array}{cc}1&0\\0&1\end{array}\right)$, and thus $M^{-1}=M^T$.}
\proof{Easily verified.}
\Lemma{\label{preserve}Transformation (\ref{trans}) preserves distances and angles between lines.}
\proof{The squared distance between any two points $U_1$ and $U_2$ is $(U_1-U_2)^T(U_1-U_2)$. The distance between the transformed points is by Lemma \ref{MTM}:
$$(\bar U_1-\bar U_2)^T(\bar U_1-\bar U_2)=(U_1-U_2)^TM^TM(U_1-U_2)=(U_1-U_2)^T(U_1-U_2)
$$	
By the cosine theorem the angles are preserved as well because the length between any two points is preserved.
}

The inverse transformation is $X=X_1+M^T\bar X$ which can be explicitly written as

\begin{equation}\label{rev}
x=a_1+\frac{a_2-a_1}{d}\bar x -\frac{b_2-b_1}{d}\bar y;~~y=b_1+\frac{a_2-a_1}{d}\bar y +\frac{b_2-b_1}{d}\bar x
\end{equation}

\subsection*{Finding the intersection points between two circles}

Let the two circles be centered at $C_1=(c_1,d_1)$ and $C_2=(c_2,d_2)$ with radii $R_1, R_2$, and the distance between the centers $d$. There are intersection points if and only if $|R_1-R_2|\le d\le R_1+R_2$.

We transform the system of coordinates by the translation and rotation (\ref{trans}) using $A_1=C_1;~A_2=C_2$.
There is actually no need to calculate the transformation. By Lemma \ref{preserve} all distances and radii are maintained.

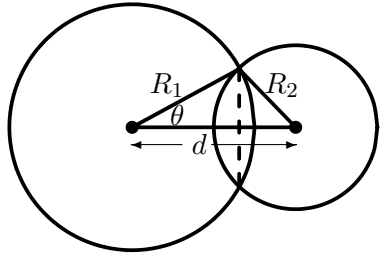
\begin{figure}[htp]
	\begin{center}\setlength{\unitlength}{1in}
		\begin{picture}(2,1.4)
		\includegraphics[width=2in]{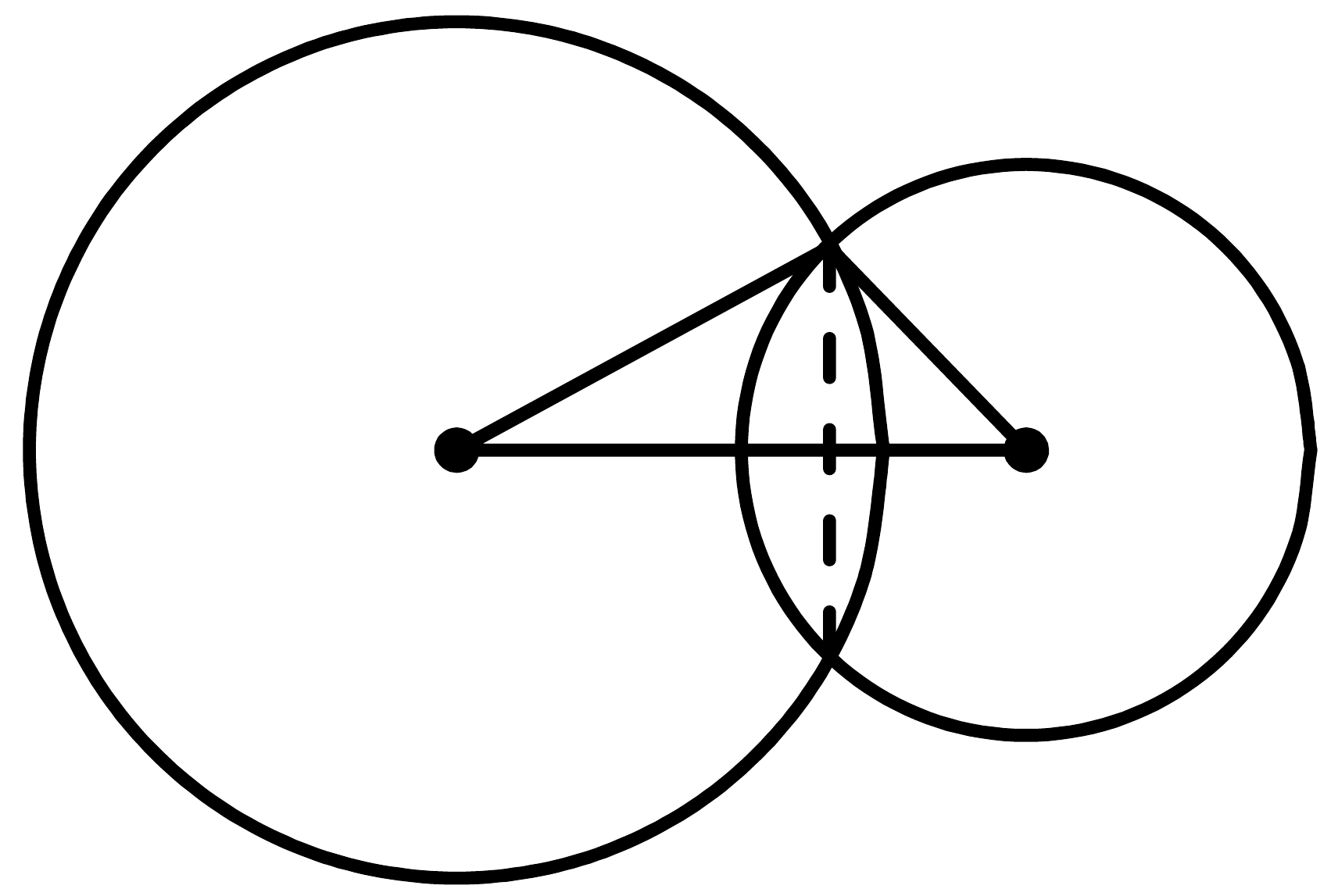}
		\put(-1.22,0.85){$R_1$}
		\put(-0.62,0.85){$R_2$}
		\put(-1.,0.53){$d$}
		\put(-1.114,0.68){$\theta$}
		\put(-1.05,0.58){\vector(-1,0){0.27}}
		\put(-0.9,0.58){\vector(1,0){0.45}}
		\end{picture}
		\caption{\label{figint} Finding the intersection points in the transformed space}
	\end{center}
\end{figure}

The two centers in the transformed space are at $(0, 0)$ and $(d, 0)$. See Figure \ref{figint}. The angle between the $x$-axis and the line connecting (0, 0) to the intersection point is defined as $\theta$. It is not required to actually calculate $\theta$. By the cosine theorem:
$$\cos\theta=\frac{R_1^2+d^2-R_2^2}{2R_1d}.$$
The intersection points in the transformed space are at
\begin{equation}\label{inters}
(R_1\cos\theta; \pm R_1\sin\theta)=\left(\frac{R_1^2+d^2-R_2^2}{2d},\pm \frac{\sqrt{[(R_1+R_2)^2-d^2][d^2-(R_1-R_2)^2]}}{2d}\right).
\end{equation}
By the condition $|R_1-R_2|\le d\le R_1+R_2$ both terms under the square root are non-negative. They cannot be both negative because $R_1+R_2>|R_1-R_2|$.  When either inequality is an equality, there is only one intersection point.
To find the intersection points in the original space, we apply (\ref{rev}) on the locations found by (\ref{inters}) using $A_1=C_1;~A_2=C_2$. The procedure requires no calculation of trigonometric functions, only one square root, when $R_1$, $R_2$, and $d$ are given.

\subsection*{Finding the intersection points between a circle and a segment}

Let the end vertices of the segment of length $d$ be at $(x_1,y_1)$ and $(x_2,y_2)$. The center of the circle is at $(x_0,y_0)$ and its radius is $R_0$.
The transformation (\ref{trans}) using $A_1=X_1;~A_2=X_2$, transfers the vertices to $(0,0)$ and $(d,0)$. The circle's center is transfered to $(\bar x_0, \bar y_0)$ and its radius is $R_0$ by Lemma \ref{preserve}. The intersection points on the $x$-axis exist if and only if $|\bar y_0|\le R_0$. They are at $\bar x=\bar x_0\pm\sqrt{R_0^2-\bar y_o^2};~~\bar y=0$. The intersection point is inside the segment, if and only if $0\le \bar x\le d$. Such a point is transformed to the original space by (\ref{rev}) using $A_1=X_1;~A_2=X_2$. It is actually simpler because $\bar y=0$.

%\newpage
 
 \subsection*{Finding whether a point is inside a triangle}
 
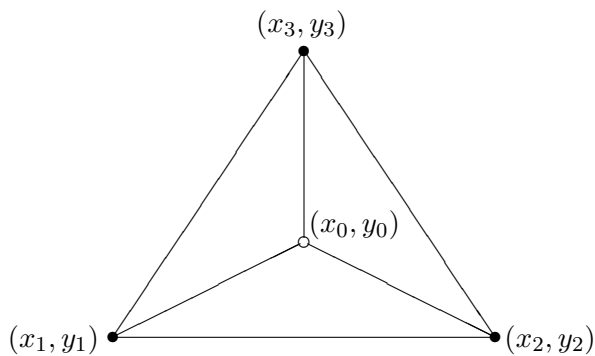
\begin{figure}[htp!]
	\begin{center}\setlength{\unitlength}{1in}
		\begin{picture}(2,2)
		\put(0,0){\circle*{0.05}}
		\put(1,1.5){\circle*{0.05}}	
		\put(2,0){\circle*{0.05}}	
		\put(1,0.5){\circle{0.05}}
\put(0,0){\line(1,0){2}}	
\put(0,0){\line(2,3){1}}	
\put(1,1.5){\line(2,-3){1}}	
\put(1,0.52){\line(0,1){0.98}}	
\put(1.02,0.49){\line(2,-1){0.98}}	
\put(0,0){\line(2,1){0.98}}
\put(-0.55,-0.05){$(x_1,y_1)$}	
\put(2.05,-0.05){$(x_2,y_2)$}
\put(0.75,1.6){$(x_3,y_3)$}	
\put(1.025,0.55){$(x_0,y_0)$}	
		\end{picture}
		\caption{\label{3T} The three triangles}
	\end{center}
\end{figure} 

\citet{pap121} suggested the following algorithm to find whether a point $(x_0, y_0)$ is inside a polygon, that may even be non-convex. The vertices of the polygon are listed in order with the last vertex identical to the first one. A side of the polygon is defined by two consecutive vertices in the correct order. For each side, the angle (between $-\pi$ and $\pi$) is defined by the direction from $(x_0, y_0)$ to the second vertex minus the direction to the first vertex. If the sum of all angles is $\pm 2\pi$ the point is inside the polygon, and if the sum is zero, it is outside.

We propose a simple procedure for a triangle $T$ which, of course, is always convex. It is based on the areas of three triangles. See Figure \ref{3T}. $(x_0,y_0)$ is inside $T$ if and only if the sum of the areas of the three triangles is equal to the area of $T$.
 There is no need for trigonometric functions that are needed for the calculation of directions. Each triangle has $(x_0,y_0)$ as one vertex and two vertices of $T$.

\Lemma{\label{area1}The area of a triangle based on three vertices $X_1, X_2, X_3$ is:
$$
\frac12\left|x_1(y_2-y_3)+x_2(y_3-y_1)+x_3(y_1-y_2)\right|
.$$}
\proof{The triangle is transformed by (\ref{trans}) using $A_1=X_1;~A_2=X_2$, so that the transformed vertices are  (0,0), $(d, 0)$ and $(\bar x_3, \bar y_3)$. The transformed triangle has the same area by Lemma \ref{preserve}. Its area is:
	$$\frac12 d|y_3|=\frac12 d\left|-\frac{y_2-y_1}{d}(x_3-x_1) +\frac{x_2-x_1}{d}(y_3-y_1)    \right|$$
which proves the lemma by collecting terms.
}

It can be verified that if the order of the points is counterclockwise, the expression inside the absolute value in Lemma \ref{area1} is positive, and it is negative if the points are in a clockwise order.

A point $(x_0,y_0)$  is inside the triangle, if and only if the sum of the three triangular areas with $(x_0,y_0)$ as a vertex and two vertices as the other two vertices is equal to the area of the triangle. Consider the four triangles in a particular order of vertices: (1, 2, 3) (the triangle $T$) and the three triangles in a specific order (0, 2, 3), (1, 0, 3) and (1, 2, 0). Suppose that the three triangles have the same sign inside the absolute value. In this case, the sum of the three areas (doubled) is the sum of three terms replacing $X_1$ by $X_0$, $X_2$ by $X_0$, etc. For example, the first term is: $x_0(y_2-y_3)+x_2(y_3-y_0)+x_3(y_0-y_2)$. If we add the three terms we get after collecting terms
$$
x_1(y_2-y_3)+x_2(y_3-y_1)+x_3(y_1-y_2)
$$
which is double the area of the triangle $T$. If the terms have different signs, the sum of the absolute values is greater than the area of the triangle.

In conclusion: if all three terms, with the agreement $X_4=X_1$, $X_5=X_2$: $$x_0(y_{i+1}-y_{i+2})+x_{i+1}(y_{i+2}-y_0)+x_{i+2}(y_0-y_{i+1})\mbox{ for }i=1,2,3$$ have the same sign, $(x_0, y_0)$ is in $T$. Otherwise, it is not.
Note that if one of the terms is zero, the point is on a side of the triangle because the area of such a triangle is zero.

	\renewcommand{\baselinestretch}{1}
	\renewcommand{\arraystretch}{1}
	\large
	\normalsize
\small

\bibliographystyle{apalike}
%\bibliography{/master}\end{document}

\begin{thebibliography}{}
	
	\bibitem[Aurenhammer and Edelsbrunner, 1984]{AE84}
	Aurenhammer, F. and Edelsbrunner, H. (1984).
	\newblock An optimal algorithm for constructing the weighted voronoi diagram in
	the plane.
	\newblock {\em Pattern Recognition}, 17:251--257.
	
	\bibitem[Aurenhammer et~al., 2013]{AKL13}
	Aurenhammer, F., Klein, R., and Lee, D.-T. (2013).
	\newblock {\em Voronoi Diagrams and Delaunay Triangulations}.
	\newblock World Scientific, New Jersey.
	
	\bibitem[Church and Garfinkel, 1978]{CG78}
	Church, R.~L. and Garfinkel, R.~S. (1978).
	\newblock Locating an obnoxious facility on a network.
	\newblock {\em Transportation Science}, 12:107--118.
	
	\bibitem[Drezner et~al., 2020]{DDK20}
	Drezner, T., Drezner, Z., and Kalczynski, P. (2020).
	\newblock Multiple obnoxious facilities location: A cooperative model.
	\newblock {\em IISE Transactions}.
	\newblock DOI: 10.1080/24725854.2020.1753898.
	
	\bibitem[Drezner, 1998]{pap121}
	Drezner, Z. (1998).
	\newblock Finding whether a point is inside a polygon and its application to
	forbidden regions.
	\newblock {\em The Journal of Management Sciences \& Regional Development},
	1:41--48.
	
	\bibitem[Drezner et~al., 2019]{DKS18}
	Drezner, Z., Kalczynski, P., and Salhi, S. (2019).
	\newblock The multiple obnoxious facilities location problem on the plane: A
	{V}oronoi based heuristic.
	\newblock {\em {OMEGA}: The International Journal of Management Science},
	87:105--116.
	
	\bibitem[Drezner and Suzuki, 2004]{DS04}
	Drezner, Z. and Suzuki, A. (2004).
	\newblock The big triangle small triangle method for the solution of non-convex
	facility location problems.
	\newblock {\em Operations Research}, 52:128--135.
	
	\bibitem[Erkut and Neuman, 1989]{Erk:89}
	Erkut, E. and Neuman, S. (1989).
	\newblock Analytical models for locating undesirable facilities.
	\newblock {\em European Journal of Operational Research}, 40:275--291.
	
	\bibitem[Gill et~al., 2005]{SNOPT}
	Gill, P.~E., Murray, W., and Saunders, M.~A. (2005).
	\newblock {SNOPT}: An {SQP} algorithm for large-scale constrained optimization.
	\newblock {\em {SIAM} {R}eview}, 47:99--131.
	
	\bibitem[Goldman and Dearing, 1975]{GD75a}
	Goldman, A.~J. and Dearing, P.~M. (1975).
	\newblock Concepts of optimal location for partially noxious facilities.
	\newblock {\em Bulletin of the Operational Research Society of America},
	23:B85.
	
	\bibitem[Hansen et~al., 1981]{HPT81}
	Hansen, P., Peeters, D., and Thisse, J.-F. (1981).
	\newblock On the location of an obnoxious facility.
	\newblock {\em Sistemi Urbani}, 3:299--317.
	
	\bibitem[Johnson, 2013]{J13}
	Johnson, R.~A. (2013).
	\newblock {\em Advanced {E}uclidean geometry}.
	\newblock Courier Corporation.
	
	\bibitem[Lee and Schachter, 1980]{LS80}
	Lee, D.~T. and Schachter, B.~J. (1980).
	\newblock Two algorithms for constructing a {D}elaunay triangulation.
	\newblock {\em International Journal of Parallel Programming}, 9:219--242.
	
	\bibitem[Okabe et~al., 2000]{OKSC00}
	Okabe, A., Boots, B., Sugihara, K., and Chiu, S.~N. (2000).
	\newblock {\em Spatial Tessellations: Concepts and Applications of {V}oronoi
		Diagrams}.
	\newblock Wiley Series in Probability and Statistics. John Wiley, Hoboken, NJ.
	
	\bibitem[Partensky, 2008]{P08}
	Partensky, M.~B. (2008).
	\newblock The circle of apollonius and its applications in introductory
	physics.
	\newblock {\em The Physics Teacher}, 46:104--108.
	
	\bibitem[Shamos and Hoey, 1975]{Sham}
	Shamos, M. and Hoey, D. (1975).
	\newblock Closest-point problems.
	\newblock {\em Proceedings 16th Annual Symposium on the Foundations of Computer
		Science, Berkeley, CA}, pages 151--162.
	
	\bibitem[Suzuki and Okabe, 1995]{SOK95}
	Suzuki, A. and Okabe, A. (1995).
	\newblock Using {V}oronoi diagrams.
	\newblock In Drezner, Z., editor, {\em Facility Location: A Survey of
		Applications and Methods}, pages 103--118. Springer, New York.
	
	\bibitem[Vorono{\"\i}, 1908]{V08}
	Vorono{\"\i}, G. (1908).
	\newblock Nouvelles applications des param{\`e}tres continus {\`a} la
	th{\'e}orie des formes quadratiques. deuxi{\`e}me m{\'e}moire. recherches sur
	les parall{\'e}llo{\`e}dres primitifs.
	\newblock {\em Journal f{\"u}r die reine und angewandte Mathematik},
	134:198--287.
	
\end{thebibliography}

\end{document}